# ON DECONVOLUTION WITH REPEATED MEASUREMENTS

By Aurore Delaigle, Peter Hall and Alexander Meister

*University of Bristol and University of Melbourne, University of Melbourne and University Stuttgart*

In a large class of statistical inverse problems it is necessary to suppose that the transformation that is inverted is known. Although, in many applications, it is unrealistic to make this assumption, the problem is often insoluble without it. However, if additional data are available, then it is possible to estimate consistently the unknown error density. Data are seldom available directly on the transformation, but repeated, or replicated, measurements increasingly are becoming available. Such data consist of "intrinsic" values that are measured several times, with errors that are generally independent. Working in this setting we treat the nonparametric deconvolution problems of density estimation with observation errors, and regression with errors in variables. We show that, even if the number of repeated measurements is quite small, it is possible for modified kernel estimators to achieve the same level of performance they would if the error distribution were known. Indeed, density and regression estimators can be constructed from replicated data so that they have the same first-order properties as conventional estimators in the known-error case, without any replication, but with sample size equal to the sum of the numbers of replicates. Practical methods for constructing estimators with these properties are suggested, involving empirical rules for smoothing-parameter choice.

**1. Introduction.** Statistical deconvolution problems arise in a great many settings, and typically have the form $g = T(f)$, where $g$ is a function about which we have data, $T$ is a transformation, and $f = T^{-1}(g)$ is a function we wish to estimate. In a large class of such problems, including density deconvolution and errors-in-variables regression, it is common to assume that $T$ is known. Indeed, the nature of the data usually precludes any other approach.

In this paper we consider cases where there is a small number replications of each intrinsically different observation, the observation errors being independent and the intrinsic parts of the observations being the same among replicates. Data of this type are numerous, and increasingly are becoming available in various fields. Examples include work of Jaech (1985), who describes an experiment where the concentration of uranium is measured for several fuel pellets; of Biemer et al. (1991), who discuss repeated observations in a social science context; of Andersen, Bro and Brockhoff (2003), on nuclear magnetic reasonance; of Bland and Altman (1986), on lung function; of Eliasziw et al. (1994), on physiotherapy for the knee; of Oman, Meir and Haim (1999), relating to kidney function; and of Dunn (1989), a brain-related study. For further medical examples, see Carroll, Ruppert and Stefanski (1995) and Dunn (2004).

When data of this type are available, it is usually possible to construct consistent estimators of the function $f$ of interest, without making parametric assumptions about the transformation $T$. We treat both density deconvolution and errors-in-variables regression, focusing on cases where the convergence rate, and first-order properties more generally, are the same when the error distribution is known and when it is not known, but is estimated from repeated measurements. In Section 2 we construct a relatively simple density estimator and generalize it to the regression case.

Theoretical properties of our estimators are taken up in Section 3. We show that a sufficient condition for first-order properties of estimators, in the cases of known and unknown error distributions, to be equivalent, is that, colloquially speaking, "the target density is smoother than half a derivative of the error density." Instances where this condition is violated are those where the convergence rate is relatively poor, even when the error density is known.

We direct attention to examples where the number of replications of each observation is relatively small. (We use the terms "replications" and "repeated measurements" synonymously.) In theoretical terms, this means that the number of replications is uniformly bounded. That is generally the case in practice, since gathering large numbers of replications is expensive in terms of time, effort or money. Moreover, particularly in cases where statistical performance is the same when the error density is known or unknown, it is seldom advantageous to have large numbers of replications.

For instance, we show that if the total number of data is $M = np$, where $p \geq 2$ equals the number of times that each of $n$ intrinsically different observations is replicated, then first-order properties of nonparametric estimators depend only on $M$, not on the separate values of $n$ and $p$. We prove this result rigorously when $p$ is bounded, but a similar argument shows that it is also valid if $p$ diverges sufficiently slowly as $M$ increases. More generally, the result holds if $M = \sum_j N_j$, where $N_j$ is the number of replicates of the $j$th intrinsically different observation. Properties of the estimator depend, to first order, only on $M$, provided that each $N_j \geq 2$.

In Section 4 we develop an adaptive, data driven procedure for smoothing-parameter choice, and show that it enjoys good performance for real and simulated datasets.

Related work in the context of density estimation includes that of Li and Vuong (1998), who derived upper bounds to convergence rates in the measurement-error problem when replications are present. Li and Vuong's results are important; they comprise some of the first contributions to density deconvolution in cases where the error distribution is not known. Nevertheless, the properties reported by Li and Vuong (1998), and bounds given also by Susko and Nadon (2002), are too coarse to permit it to be shown that convergence rates can be identical in the cases of known and unknown error distributions. Further discussion is given in Section 3.5.

Recent, related research in the regression setting, and in the econometrics literature, includes that of Li (2002), Li and Hsiao (2004) and Schennach (2004a, 2004b), who demonstrated that replications can be used to good effect in regression problems with measurement error. See also the work of Horowitz and Markatou (1996) on error estimation from panel data, and the extensive literature, accessible through the work of Newey and Powell (2003), on inference in the context of instrumental variables. However, except in parametric contexts, this and related work is not sufficiently detailed to show that the convergence rates familiar in problems where the error distribution is known can also be enjoyed when the distribution is accessible only via repeated measurements.

The problem of density estimation with unknown error density, estimated from a sample of the error, has been considered by Diggle and Hall (1993), Barry and Diggle (1995) and Neumann (1997). Madansky (1959), Carroll, Eltinge and Ruppert (1993) and Huang and Yang (2000), among others, have discussed linear regression with replicated data, when at least some of the predictors are measured with error. Early work on the problem of density deconvolution, under the assumption of known distribution of measurement error, includes that of Carroll and Hall (1988), Stefanski and Carroll (1990) and Fan (1991). More recent contributions, including surveys of earlier research, include the papers of Delaigle and Gijbels (2002, 2004) and van Es and Uh (2005). The literature on kernel methods for errors-in-variables regression is particularly large, and is surveyed by Carroll, Ruppert and Stefanski (1995).

## 2. Models and methodology.

2.1. *Density deconvolution.* Suppose we observe

(2.1) $\qquad W_{jk} = X_j + U_{jk} \qquad \text{for } 1 \leq k \leq N_j \text{ and } 1 \leq j \leq n,$

where the random variables $X_j$ are identically distributed as $X$, the $U_{jk}$'s are identically distributed as $U$, and the $X_j$'s and $U_{jk}$'s are totally independent. We wish to estimate the density of $X$. In the context of our discussion in Section 1, (2.1) indicates that there are $n$ subsets of "intrinsically different" data and, within the $j$th of these subsets, $N_j$ repeated, or replicated, measurements of the variable $X_j$.

Let $f_U$ and $f_X$ denote the respective densities of $U$ and $X$, and write $f_U^{\text{Ft}}$ and $f_X^{\text{Ft}}$ for the respective characteristic functions (i.e., the Fourier transforms of those densities). Provided that

(2.2) $f_U^{\text{Ft}}$ is real-valued and does not vanish at any point on the real line,

a consistent estimator of $f_U^{\text{Ft}}$ is given by

$$（2.3）\quad \hat{f}_U^{\text{Ft}}(t) = \left| \frac{1}{N} \sum_{j=1}^n \sum_{(k_1,k_2) \in \mathcal{S}_j} \cos\{t(W_{jk_1} - W_{jk_2})\} \right|^{1/2},$$

where $\mathcal{S}_j$ denotes the set of $\frac{1}{2} N_j(N_j - 1)$ distinct pairs $(k_1, k_2)$ with $1 \leq k_1 < k_2 \leq N_j$, $N = N(n) = \frac{1}{2} \sum_{j \leq n} N_j(N_j - 1)$, and we ignore values of $j$ for which $N_j = 1$. Assumption (2.2) is conventional when using kernel methods for density deconvolution; see Stefanski and Carroll (1990) and Fan (1991), for example.

An estimator of $f_X$ is given by

$$\hat{f}_X(x) = \frac{1}{Mh} \sum_{j=1}^n w_j \sum_{k=1}^{N_j} \widehat{L}\left( \frac{x - W_{jk}}{h} \right),$$

where $M = \sum_j N_j$, the weights $w_j$ are nonnegative and satisfy $\sum_j w_j N_j = M$,

$$（2.4）\quad \widehat{L}(u) = \frac{1}{2\pi} \int e^{-itu} \frac{K^{\text{Ft}}(t)}{\hat{f}_U^{\text{Ft}}(t/h) + \rho} \, dt,$$

$K$ is a symmetric kernel function with compactly supported Fourier transform $K^{\text{Ft}}$, $h > 0$ is a bandwidth, and $\rho \geq 0$ is a ridge parameter.

We introduce the ridge only so we can take expectation without concern for fluctuations of the denominator in the integral at (2.4). The ridge would not be necessary if our aim were to develop limit theory for $\hat{f}_X$ that did not involve taking expected values. See Section 3.1 for discussion and theory in the case $\rho = 0$.

If $f_U$ were known then, instead of $\hat{f}_X$, we would use the following generalization of the conventional deconvolution estimator:

$$\tilde{f}_X(x) = \frac{1}{Mh} \sum_{j=1}^n w_j \sum_{k=1}^{N_j} L\left( \frac{x - W_{jk}}{h} \right)$$

[see, e.g., Carroll and Hall (1988)], where

$$L(u) = \frac{1}{2\pi} \int e^{-itu} \frac{K^{\mathrm{Ft}}(t)}{f_U^{\mathrm{Ft}}(t/h)} \, dt.$$

The bias of $\tilde{f}_X$ does not depend on choice of the weights, and it can readily be shown that the asymptotic variance is minimized by taking each $w_j = 1$. Optimality of this choice persists in the case of regression deconvolution, which we consider in Section 2.2.

Therefore, we take each $w_j = 1$ in the work below. In particular, $\hat{f}_X$ and $\tilde{f}_X$ henceforth denote the estimators

$$\hat{f}_X(x) = \frac{1}{Mh} \sum_{j=1}^{n} \sum_{k=1}^{N_j} \hat{L}\left(\frac{x - W_{jk}}{h}\right)$$

and

$$\tilde{f}_X(x) = \frac{1}{Mh} \sum_{j=1}^{n} \sum_{k=1}^{N_j} L\left(\frac{x - W_{jk}}{h}\right).$$

Section 3.3 demonstrates that $\hat{f}_X$ is first-order equivalent to $\tilde{f}_X$. For this result, and in the setting of "ordinary-smooth errors" [see (3.1)], the main assumption needed is that $f_X$ be sufficiently smooth relative to $f_U$. See condition (3.12). Properties of $\tilde{f}_X$ are summarized in Section 3.4.

2.2. *Errors-in-variables regression.* Here the model at (2.1) is extended, so that it addresses data $(W_{jk}, Y_j)$ generated as

(2.5) $$W_{jk} = X_j + U_{jk}, \qquad Y_j = g(X_j) + V_j,$$

$$\text{for } 1 \leq k \leq N_j \text{ and } 1 \leq j \leq n,$$

where the $X_j$'s, $U_{jk}$'s and $V_j$'s are identically distributed as $X$, $U$ and $V$, respectively, $E(V) = 0$, $E(V^2) < \infty$, and the $X_j$'s, $U_{jk}$'s and $V_j$'s are totally independent. We wish to estimate the function $g$.

Define

(2.6)
$$\hat{a}(x) = \frac{1}{Mh} \sum_{j=1}^{n} \sum_{k=1}^{N_j} Y_j \hat{L}\left(\frac{x - W_{jk}}{h}\right),$$

$$\tilde{a}(x) = \frac{1}{Mh} \sum_{j=1}^{n} \sum_{k=1}^{N_j} Y_j L\left(\frac{x - W_{jk}}{h}\right).$$

In the classical case, where $f_U$ is known and each $N_j = 1$, the standard kernel estimator of $g$ is $\tilde{g} = \tilde{a}/\tilde{f}_X$ and, of course, $\tilde{g}$ is also appropriate in the case of replicated data.

The intuition behind $\tilde{g}$ is that $\tilde{a}$ is a consistent estimator of the function $a = f_X g$. When $f_U$ is not known we can estimate $a$ by $\hat{a}$, and so we can modify $\tilde{g}$ in the manner of Section 2.1, estimating $g$ by $\hat{g} = \hat{a}/\hat{f}_X$. We show in Section 3.6 that $\hat{g}$ is first-order equivalent to $\tilde{g}$.

## 3. Theoretical properties.

3.1. *Density deconvolution.* First we state assumptions. We ask that, for constants $\alpha > 0$ and $B_1 > 1$, and all real $t$,

(3.1) $$B_1^{-1}(1+|t|)^{-\alpha} \leq |f_U^{\mathrm{Ft}}(t)| \leq B_1(1+|t|)^{-\alpha}.$$

This is often referred to as the case of ordinary-smooth errors. The importance of the lower bound in (3.1), in addition to the upper bound (which is conventional when deriving convergence rates), is discussed in Section 3.3.

Given $\beta, B_2 > 0$, let $\mathcal{F}(\beta, B_2)$ denote the class of densities $f_X$ for which

$$\sup_{-\infty < t < \infty} (1+|t|)^\beta |f_X^{\mathrm{Ft}}(t)| \leq B_2.$$

[The class $\mathcal{F}(\beta, B_2)$ is a Fourier analogue of Fan's class $\mathcal{C}_{m,\alpha,B}$ of functions; his $m + \alpha + 1$ is our $\beta$.] Let $K$ have the property

(3.2) $\sup |K^{\mathrm{Ft}}| < \infty$ and, for some $c > 0$, $K^{\mathrm{Ft}}(t) = 0$ for all $|t| > c$.

The kernels used in deconvolution commonly have this property, and so, while our results can be derived under weaker conditions, there is little motivation for that generalization.

The theorem below gives an upper bound to pointwise mean-squared distance between $\hat{f}_X$ and $\tilde{f}_X$, uniformly in all points and all densities $f_X \in \mathcal{F}(\beta, C_2)$. In Section 3.3 we use that result to show that, if the bandwidth $h$ is chosen so that it gives optimal performance of $\hat{f}_X$, and if a relation (3.12) on the relative smoothnesses of $f_U$ and $f_X$ holds, then the difference between $\hat{f}_X$ and $\tilde{f}_X$ is negligible relative to the distance between either estimator and the true density, $f_X$.

THEOREM 3.1. *Let $C_1 > 1$ and $C_2, \beta > 0$. Assume that* (i) $1 \leq N_j \leq C_1$ *for each $j$;* (ii) $N(n) \geq C_1^{-1} n$ *for each $n \geq 1$;* (iii) $f_U^{\mathrm{Ft}}$ *satisfies (3.1);* (iv) $\alpha > \frac{1}{2}$; (v) $K^{\mathrm{Ft}}$ *satisfies (3.2);* (vi) $h_1(n) \leq h \leq h_2(n)$, *where $h_2(n) \to 0$ and, for some $\delta > 0$, $n^{(1-\delta)/4\alpha} h_1(n)$ is bounded away from zero; and* (vii) $c_1 n^{-c_2} \leq \rho \leq c_3 \min\{h_1(n)^{4\alpha+2}, n^{-1}\}$, *where $c_1, c_2, c_3 > 0$. Then, for each integer $k \geq 1$,*

(3.3) $$\sup_{f_X \in \mathcal{F}(\beta, C_2)} \sup_{-\infty < x < \infty} E\{\hat{f}_X(x) - \tilde{f}_X(x)\}^2 \leq \mathrm{const.} p_n,$$

*where, for each integer $k \geq 1$,*

$$
\begin{aligned}
(3.4) \quad p_n = p_n(k) = {} & n^{-1}\{h^{\beta-2\alpha-1} + h^{2(\beta-2\alpha)-1} + (\log n)^2\} \\
& + n^{-2}(h^{2(\beta-4\alpha)-2} + h^{-6\alpha-1}) + n^{-k}h^{-4(k+2)\alpha-2}
\end{aligned}
$$

*and the constant in (3.3) depends on $k$ but not on $h \in [h_1(n), h_2(n)]$ or on $n$.*

Technical arguments are given in a longer version of this paper [Delaigle, Hall and Meister (2006)]. Theorem 3.1 remains correct without condition (i), that is, without the assumption that the $N_j$'s are bounded uniformly in $j$ and $n$. However, if (i) is dropped, then the asymptotic properties of $\hat{f}$ cannot be discussed simply in terms of the size of $M$, and that difficulty hampers elucidation of our results. Indeed, if condition (i) is removed, then, depending on the size of the $N_j$'s, and on the frequency with which large $N_j$'s occur, properties of $\hat{f}$ can be very close to those of a standard kernel estimator based on the (unobservable) data $X_j$. In practice, the expense, in terms of time, effort or money, of making repeated measurements usually ensures that the $N_j$'s are relatively small, typically no more than 2 to 5, and so we shall retain condition (i).

We argued in Section 2 that, if we were to develop a limit theory that did not involve taking expected values, the ridge parameter $\rho$ could be taken equal to zero. In that setting we should replace uniform pointwise error, at (3.3), by error at a single point, or by a global metric such as integrated squared error. Otherwise, we incur a logarithmic penalty on the right-hand side of (3.3). [This is to be expected, since the same penalty arises in more conventional problems; see, e.g., Bickel and Rosenblatt (1973).] We should also remove the supremum over densities $f_X \in \mathcal{F}(\beta, C_2)$, since the uniformity implied by the supremum is not meaningful if we remove the expectation.

For the sake of definiteness, when working with $\rho = 0$, we measure accuracy in terms of squared error at a particular point, or integrated squared error. To treat the latter, note that (3.3) implies that, for each pair $x_1, x_2$ for which $-\infty < x_1 < x_2 < \infty$,

$$
(3.5) \quad \sup_{f_X \in \mathcal{F}(\beta, C_2)} \int_{x_1}^{x_2} E\{\hat{f}_X(x) - \tilde{f}_X(x)\}^2 \, dx = O(p_n).
$$

Let $\hat{f}_X^0(x)$ denote the version of $\hat{f}_X$ constructed with $\rho = 0$. We claim that (3.5) continues to apply to $\hat{f}_X^0$, provided the expectation and supremum over $f_X$ are removed from the left-hand side, and the right-hand side is interpreted in an "in probability" sense. Moreover, squared error at each fixed point $x$ converges at the same rate:

$$
|\hat{f}_X^0(x) - \tilde{f}_X(x)| = O_p(p_n^{1/2}),
$$

(3.6)
$$\int_{x_1}^{x_2} \{\hat{f}_X^0(x) - \tilde{f}_X(x)\}^2\, dx = O_p(p_n).$$

THEOREM 3.2. *Let $C_1 > 1$, let $C_2, \alpha, \beta > 0$, let $-\infty < x_1 < x_2 < \infty$, and take $\rho = 0$ in the definition of $\widehat{L}$, at (2.4), and hence, also in the definition of $\hat{f}_X$, obtaining the estimator $\hat{f}_X^0$. Assume that conditions* (i)–(vi) *in Theorem 3.1 hold. Then (3.6) holds for each $f_X \in \mathcal{F}(\beta, C_2)$, each $x \in (-\infty, \infty)$ and each pair $x_1, x_2$ for which $-\infty < x_1 < x_2 < \infty$.*

Results (3.6) and (3.10), below, show that optimal convergence rates can be achieved using a single smoothing parameter, the bandwidth, rather than two parameters, the bandwidth and ridge.

3.2. *Asymptotic optimality.* The size of bandwidth that minimizes pointwise mean squared error, when using $\tilde{f}_X$ to estimate $f_X$, is $h \asymp h_0 \equiv n^{-1/\{2(\alpha+\beta)-1\}}$; and, for such a bandwidth, pointwise mean squared error of $\tilde{f}_X$ is of size $q_n$, where

(3.7)  $$q_n = n^{-2(\beta-1)/\{2(\alpha+\beta)-1\}}.$$

The same result holds if we replace $\tilde{f}_X$ by the errors-in-variables regression estimator, $\tilde{g}$, which we define in Section 3.6. See Fan (1991) and Fan and Truong (1993) for discussion of theory in these respective cases, and also for proofs of lower bounds which show that the rate $q_n$ is minimax optimal, in an $L_2$ sense.

However, these results address only the case where there is no replication, that is, each $N_j = 1$. In the case of upper bounds, generalization to settings where each $N_j \geq 2$ is relatively straightforward. See Section 3.4 for details. Below we generalize lower bounds in the setting of density deconvolution.

THEOREM 3.3. *Assume that $\alpha, \beta > \frac{1}{2}$. Let $\mathcal{F}(\beta, C)$ denote the class of densities $f_X$ defined in Section 3.1, and write $\breve{\mathcal{F}}$ for the class of all measurable functionals of the data. Assume that $2 \leq N_j \leq B$ for each $j$, where $2 \leq B < \infty$. Then, for each fixed $x$ and each sufficiently large $C > 0$, there exists $D > 0$ such that, for all sufficiently large $n$,*

(3.8)  $$\inf_{\breve{f} \in \breve{\mathcal{F}}} \sup_{f_X \in \mathcal{F}(\beta, C)} E_{f_X}\{\breve{f}(x) - f_X(x)\}^2 \geq D q_n.$$

3.3. *Equivalence of $\hat{f}_X$ and $\tilde{f}_X$.* In view of the results given in Section 3.2, and in order to establish that $\hat{f}_X$ is asymptotically equivalent to $\tilde{f}_X$ when the latter is performing optimally, it is instructive to show that when $h \asymp h_0$,

(3.9)  $$\sup_{f_X \in \mathcal{F}(\beta, C_2)} \sup_{-\infty < x < \infty} E\{\hat{f}_X(x) - \tilde{f}_X(x)\}^2 = o(q_n),$$

if the ridge-prameter $\rho$ is taken to be nonzero; or, if the ridge is zero, that

$$|\hat{f}_X^0(x) - \tilde{f}_X(x)| = o_p(q_n^{1/2}),$$

(3.10)
$$\int_{x_1}^{x_2} \{\hat{f}_X^0(x) - \tilde{f}_X(x)\}^2 \, dx = o_p(q_n).$$

Compare with (3.6). In fact, (3.9) and (3.10) follow from Theorems 3.1 and 3.2, respectively, if we prove that

(3.11) $$q_n = o(p_n).$$

Provided

(3.12) $$\beta > \alpha + \tfrac{1}{2},$$

it is straightforward to show that if $h \asymp h_0$, then

(3.13)
$$n^{-1}\{h^{\beta-2\alpha-1} + h^{2(\beta-2\alpha)-1} + (\log n)^2\}$$
$$+ n^{-2}(h^{2(\beta-4\alpha)-2} + h^{-6\alpha-1}) = o(q_n),$$

and also that if $k$ is sufficiently large and $h \asymp h_0$, then $n^{-k}h^{-4(k+2)\alpha-2} = o(q_n)$. This result and (3.13) imply (3.11).

Therefore, condition (3.12), which can be characterized colloquially as the assertion that "$f_X$ is smoother than half a derivative of $f_U$," is sufficient to ensure that, in deconvolution problems, there is no first-order loss of performance in using replicated data to estimate the error density when the latter is not known. Intuition behind (3.12) is given in Section 3.5.

Of course, (3.12) fails if $\alpha$ is too large; that is, if $f_U$ is too smooth. This is the reason for placing the lower bound on $|f_U^{\mathrm{Ft}}(t)|$ in (3.1). Without that bound, $f_U$ can be arbitrarily smooth. It can be shown that if $\beta < \alpha$, then $\hat{f}_X$ and $\tilde{f}_X$ are not asymptotically equivalent, and the minimax-optimal, pointwise convergence rate of an estimator of $f_X$ can be no faster than $n^{-2(\beta-1)/(4\alpha-1)}$, which is strictly slower than the rate of convergence of $\tilde{f}_X$ to $f_X$. However, the case where $\alpha \leq \beta \leq \alpha + \tfrac{1}{2}$ is still unclear.

3.4. *Properties of $\tilde{f}_X$.* Let $\check{f}_X$ denote the "standard" version of $\tilde{f}_X$, obtained by taking $N_j = 1$ for each $j$, but with sample size $M$ rather than $n$. Theorem 3.4, which is given below and is straightforward to derive, argues that the bias of $\tilde{f}_X$ is identical to that of $\check{f}_X$, and that the variance of $\tilde{f}_X$ equals that of $\check{f}_X$, to first order.

Recall that $U$ and $X$ have the distributions of $U_{jk}$ and $X_j$, respectively, that $W = X + U$, and that $N = \tfrac{1}{2}\sum_{j \leq n} N_j(N_j - 1)$. Put

$$m_n(x) = \int K(u) f_X(x - hu) \, du,$$

$$v_n(x) = \frac{1}{M}\left\{\frac{1}{h}\int L(u)^2 f_W(x-hu)\,du - m_n(x)^2\right\},$$

$$w_n(x) = \frac{2N}{M^2}\left\{\frac{1}{h}\int K(u)^2 f_X(x-hu)\,du - m_n(x)^2\right\}.$$

THEOREM 3.4. *The mean and variance of $\check{f}_X(x)$ equal $m_n(x)$ and $v_n(x)$, respectively; the mean of $\tilde{f}_X(x)$ equals $m_n(x)$; and the variance of $\tilde{f}_X(x)$ equals $v_n(x) + w_n(x)$.*

The quantity $w_n$ is generally of strictly smaller order than $v_n$, since $\int K^2$ remains fixed but $\int L^2$ diverges as $h$ decreases. Therefore, in terms of first-order properties of mean and variance, $\tilde{f}_X$ and $\check{f}_X$ have identical performance. In view of this property, and bearing in mind the asymptotic equivalence of $\hat{f}_X$ and $\tilde{f}_X$ noted in Section 3.3, we can fairly say that:

(3.14) to first order, $\hat{f}_X$ has the same properties as a conventional deconvolution density estimator, computed when the error density is known and the sample size is $M$ but without any replication.

Of course, this assertion requires (3.9) and, hence, needs (3.12).

Together, (3.8), (3.9) and (3.14) demonstrate minimax optimality of the estimator $\hat{f}_X$. Of course, this property necessitates the supremum being taken over $f_X$ in (3.9). That requirement motivated our introduction of the ridge parameter in our definition of $\hat{f}_X$.

3.5. *Discussion of different approaches to density deconvolution.* Let $(2.2)'$ denote the version of (2.2) where the assumption that $f_U^{\mathrm{Ft}}$ is real-valued is omitted. For cases where $(2.2)'$ holds but (2.2) fails, Li and Vuong (1998) suggest an estimator of $f_U^{\mathrm{Ft}}$ quite different from our $\hat{f}_U^{\mathrm{Ft}}$. However, from a practical viewpoint, the condition that $f_U^{\mathrm{Ft}}$ be real-valued is mild. In particular, in the nonparametric literature on density deconvolution and errors-in-variables regression where $f_U$ is assumed known, that quantity is invariably taken to be symmetric, in which case $f_U^{\mathrm{Ft}}$ is real-valued.

The alternative estimator suggested by Li and Vuong (1998) in the context of $(2.2)'$ requires the distributions of both $U$ and $X$ to have characteristic functions that do not vanish anywhere (see Li and Vuong's condition A3) and also to be compactly supported (see their assumption A4). We are not aware of a distribution which enjoys both these properties. Certainly, none of the standard, compactly-supported distributions satisfy A3. This, and the numerical complexity of Li and Vuong's estimator, discouraged us from considering their technique.

If $\alpha$ is sufficiently less than $\beta$, then the problem of estimating $f_U$ from the differences $W_{jk_1} - W_{jk_2}$ is more difficult statistically, although more straightforward numerically, than the problem of estimating $f_U$ from the raw data $W_{jk}$. This indicates why condition (3.12) is required. For values of $\alpha$ that are large relative to $\beta$, alternative deconvolution methods may possibly give better theoretical performance, although we are not aware of any that are attractive computationally.

3.6. *Errors-in-variables regression.* The results in this section are closely analogous to those in earlier sections, so we give only an outline. Recall from Section 2.2 that, under the model (2.5), our estimator of $g$ is $\hat{g} = \hat{a}/\hat{f}_X$, where $\hat{a}$ is an estimator, defined at (2.6), of $a = f_X g$. Properties of $\hat{g}$ follow directly from those of the numerator and denominator in the ratio $\hat{a}/\hat{f}_X$. The denominator is treated in Theorems 3.1 and 3.2; here we address the numerator.

Given $f_X \in \mathcal{F}(\beta, C_2)$, let $\mathcal{G}(\beta, C_2 | f_X)$ denote the class of functions $g$ for which

$$\sup_{-\infty < t < \infty} (1 + |t|)^\beta \left| \int e^{itx} f_X(x) g(x) \, dx \right| \leq C_2.$$

Recall that conditions associated with the errors-in-variables model (2.5) include the assumption that $E(V) = 0$ and $E(V^2) < \infty$.

THEOREM 3.5. *Let $C_1 > 1$ and $C_2, \alpha, \beta > 0$. Assume* (i)–(vii) *in Theorem 3.1. Then, for each integer $k \geq 1$,*

$$(3.15) \quad \sup_{f_X \in \mathcal{F}(\beta,C_2), g \in \mathcal{G}(\beta,C_2|f_X)} \sup_{-\infty < x < \infty} E\{\hat{a}(x) - \tilde{a}(x)\}^2 \leq \mathrm{const.} p_n,$$

*where $p_n$ is as at (3.4) and the constant in (3.15) depends on $k$ but not on $h \in [h_1(n), h_2(n)]$ or on $n$.*

We know from Section 3.3 that, if $\alpha$ and $\beta$ satisfy (3.12), and if $h$ is of the same size as the bandwidth that minimizes mean squared error of $\tilde{f}_X$ (this is also the size of the optimal bandwidth for $\tilde{a}$ and $\tilde{g}$), then $p_n = o(q_n)$. [Recall that $q_n$ is given by (3.7), and that $q_n^{1/2}$ equals the minimum order of magnitude of error for estimators of $f_X$, $a$ and $g$.] It then follows from Theorems 3.1 and 3.5, and (3.11), that if conditions (i)–(vii) hold, $\hat{f}_X(x) - \tilde{f}_X(x) = o_p(q_n^{1/2})$ and $\hat{a}(x) - \tilde{a}(x) = o_p(q_n^{1/2})$. Therefore, provided $f_X(x) > 0$, we have

$$(3.16) \qquad \hat{g}(x) = \frac{\hat{a}(x)}{\hat{f}_X(x)} = \frac{\tilde{a}(x)}{\tilde{f}_X(x)} + o_p(q_n^{1/2}) = \tilde{g}(x) + o_p(q_n^{1/2}).$$

That is, if the bandwidth is chosen so that it is optimal for estimating $g$ by $\tilde{g}$, then $\hat{g}$ is first-order equivalent to $\tilde{g}$.

It is straightforward to state and prove the analogue of Theorem 3.4 for the estimator $\tilde{a}$ instead of $\tilde{f}_X$. This leads directly to the analogue of (3.14), where the only change necessary is to replace $\hat{f}_X$ by $\hat{g}_X$ and alter "density estimator" to "regression estimator."

An argument similar to that used in Section 5 to derive Theorem 3.2 can be employed to show that (3.16) holds even if the ridge parameter, $\rho$, is taken as zero. Therefore, (3.14) applies in the ridge-free case.

3.7. *Supersmooth error case.* All our discussion in the previous paragraphs was based on the assumption that the error distribution is ordinary smooth, and, in particular, satisfies (3.1). It is also of interest to treat the case of supersmooth errors, so named because there the error density is infinitely differentiable. In that context the following condition is imposed in place of (3.1): for constants $\alpha > 0$, $\gamma > 0$ and $B_1 > 1$, and all real $t$,

$$(3.17) \qquad B_1^{-1} \exp(-\gamma |t|^\alpha) \leq |f_U^{\mathrm{Ft}}(t)| \leq B_1 \exp(-\gamma |t|^\alpha).$$

For such error distributions, pointwise mean squared error, when employing $\tilde{f}_X$ to estimate $f_X$, is of optimal order when using a bandwidth $h = D(\log n)^{-1/\alpha}$, where $D > (4\gamma)^{1/\alpha}$ denotes a constant. In this case, pointwise mean squared error of $\tilde{f}_X$ is of size $q_n = (\log n)^{-2(\beta-1)/\alpha}$. Here, the rate of convergence of the estimator $\tilde{f}_X$ is so slow that the loss of performance incurred by estimating $f_U$ from the data, and using $\hat{f}_X$ instead of $\tilde{f}_X$, is negligible, regardless of restrictions such as (3.12). In particular, the following theorem holds. Its proof follows the lines of that of Theorem 3.1, but is more straightforward.

THEOREM 3.6. *Let $C_1 > 1$ and $C_2, C_3, \alpha, \beta, \gamma > 0$. Assume that* (i) $1 \leq N_j \leq C_1$ *for each $j$;* (ii) $N(n) \geq C_1^{-1} n$ *for each $n \geq 1$;* (iii) $f_U^{\mathrm{Ft}}$ *satisfies (3.17);* (iv) $K^{\mathrm{Ft}}$ *satisfies (3.2) with $c = 1$;* (v) $h = D(\log n)^{-1/\alpha}$, *with $D > (4\gamma)^{1/\alpha}$; and* (vi) $\rho = C_2 n^{-\kappa}$, *with $\kappa > \frac{1}{4}$. Then, for some $\epsilon > 0$,*

$$\sup_{f_X \in \mathcal{F}(\beta, C_3)} \sup_{-\infty < x < \infty} E\{\hat{f}_X(x) - \tilde{f}_X(x)\}^2 \leq \mathrm{const.} n^{-\epsilon}.$$

This result is readily generalized to the estimator $\hat{g}$, provided $h$ is chosen so that the optimal convergence rate for $\tilde{g}$ as an estimator of $g$ is attained. In particular, if $h = D(\log n)^{-1/\alpha}$ where $D > (4\gamma)^{1/\alpha}$, then $\hat{g}$ is first-order equivalent to $\tilde{g}$.

## 4. Numerical properties.

4.1. *Simulated examples.* We study numerical properties of the estimators $\hat{f}_X$ and $\hat{g}$ in several simulated examples. In the density case, and follow-ing model (2.1), we generate 500 random samples of replicated observations for $n$ individuals, $W_{ij}$, where $i = 1, \ldots, n$ and $j = 1, \ldots, N_i$. We take the noise-to-signal ratio $\sigma_U^2/\sigma_X^2$ equal to 25%, except in the case of density (iii) below, where we take $\sigma_U^2/\sigma_X^2 = 10\%$. The notation $\sigma_T^2$ denotes the variance of a random variable $T$. The error density $f_U$ is chosen to be a Laplace or a centered normal density. In each instance where the first of these choices is used, (3.12) is satisfied; the second choice corresponds to a supersmooth density, and there (3.12) is not relevant.

We consider four target densities $f_X$: (i) $X \sim 0.5\text{N}(-3,1) + 0.5\text{N}(2,1)$, (ii) $X \sim \chi^2(3)$, (iii) $X \sim \sum_{\ell=0}^{5}(2^{5-\ell}/63)\text{N}\{65 - 96 \times 2^{-\ell}/21, (32/63)^2/2^{2\ell}\}$ and (iv) $X \sim \text{N}(0,1)$. Density (i) is bimodal and symmetric, density (ii) is asymmetric and density (iii) is the smooth comb density discussed by Marron and Wand (1992). Note that, even in the error-free case, the latter density is particularly hard to estimate because of its numerous features.

In the regression case we generate 500 datasets of randomly-sampled vectors $(W_{ij}, Y_i)$, $i = 1, \ldots, n$, $j = 1, \ldots, N_i$, according to the model (2.5). The density $f_X$ is chosen to be a uniform $U[0,1]$ or a normal $\text{N}(0.5, \sigma_X^2)$ density, with $\sigma_X^2$ chosen so that 0 and 1 are respectively the 0.025 and 0.975 quantiles of $f_X$. The error density $f_U$ is a Laplace or centered normal density, and the noise-to-signal ratio $\sigma_U^2/\sigma_X^2$ equals 10%. Except for our Bernoulli regression example [see case (iii) below], the error density $f_V$ is a centered normal density such that the noise-to-signal ratio $\sigma_V^2/\sigma^2(g)$ equals 10%, where $\sigma^2(g)$ denotes the mean squared deviation of $g$ from its average value.

We consider three regression curves: (i) $g(x) = x^2(1-x)^2$, (ii) $g(x) = 3x + 20(2\pi)^{-1/2}\exp\{-100(x-\frac{1}{2})^2\}$, (iii) $Y|X = x \sim \text{Bernoulli}\{g(x)\}$, with $g(x) = 0.45\sin(2\pi x) + 0.5$. Note that curve (i) is unimodal and symmetric around 0.5, curve (ii) is a mixture of a straight line and an exponential curve, and curve (iii) is an asymmetric sinusoid.

We sought an automatic way of choosing the bandwidth, $h$. In the density case, we suggest using $\hat{h}_{\text{PI}}$, the plug-in bandwidth of Delaigle and Gijbels (2002, 2004), where the characteristic function of the error is replaced by (2.3). This procedure is justified by the discussion in Section 3.3. In the regression case, a bandwidth-choice procedure could also be based on a data-driven selector for the known error case. However, since, to our knowledge, there does not exist such a method, we must first propose one.

A cross-validation (CV) criterion for selecting $h$ would choose

$$h_{\text{CV}} = \arg\min_h \sum_{k=1}^{n} \left( \frac{Y_k - \sum_{j=1}^{n} Y_j S_j(X_k)}{1 - S_k(X_k)} \right)^2,$$

where, for $j = 1, \ldots, n$,

$$S_j(x) = \sum_{\ell=1}^{N_j} L\left(\frac{x - W_{j\ell}}{h}\right) \bigg/ \sum_{J=1}^{n} \sum_{\ell=1}^{N_J} L\left(\frac{x - W_{J\ell}}{h}\right).$$

Since the observations $X_k$ are not available, we need to replace all quantities of the form

$$L\left(\frac{X_k - W_{j\ell}}{h}\right) = \frac{1}{2\pi} \int \exp(-itX_k/h) \exp(itW_{j\ell}/h) \frac{K^{\mathrm{Ft}}(t)}{f_U^{\mathrm{Ft}}(t/h)} \, dt,$$

by empirical estimators. We suggest replacing $\exp(-itX_k/h)$ by an estimator of its expected value, $f_X^{\mathrm{Ft}}(-t/h)$, based on the replications of the $k$th intrinsic observation. Such an estimator can be defined by $\hat{f}_W^{\mathrm{Ft}}(-t/h)/f_U^{\mathrm{Ft}}(-t/h)$, where $\hat{f}_W^{\mathrm{Ft}}(t) = K^{\mathrm{Ft}}(ht) \sum_{m=1}^{N_k} \exp(itW_{km})$ is a kernel estimator of $f_W^{\mathrm{Ft}}$. Proceeding that way, our CV criterion becomes

$$\text{(4.1)} \qquad \tilde{h}_{\mathrm{CV}} = \arg\min_h \sum_{k=1}^n \left( \frac{Y_k - \sum_{j=1}^n Y_j \widehat{S_j(X_k)}}{1 - \widehat{S_k(X_k)}} \right)^2,$$

where

$$\text{(4.2)} \quad \widehat{S_j(X_k)} = \sum_{m=1}^{N_k} \sum_{\ell=1}^{N_j} L_2\left(\frac{W_{km} - W_{j\ell}}{h}\right) \Big/ \sum_{J=1}^n \sum_{m=1}^{N_k} \sum_{\ell=1}^{N_J} L_2\left(\frac{W_{km} - W_{J\ell}}{h}\right),$$

with $L_2(x) = (2\pi)^{-1} \int \exp(-itx/h) |K^{\mathrm{Ft}}(t)|^2 |f_U^{\mathrm{Ft}}(t/h)|^{-2} \, dt$.

In the case of unknown error density, we define $\hat{h}_{\mathrm{CV}}$ as in (4.1) but we replace $L_2$ in (4.2) by

$$\widehat{L}_2(x) = (2\pi)^{-1} \int \exp(-itx/h) |K^{\mathrm{Ft}}(t)|^2 |\hat{f}_U^{\mathrm{Ft}}(t/h)|^{-2} \, dt,$$

with $\hat{f}_U^{\mathrm{Ft}}(t)$ as in (2.3). As in the error-free case, the computations needed to calculate this bandwidth can be reduced considerably by binning the data. See, for example, Fan and Gijbels (1996), page 96. We suggest placing the $W_{ij}$'s into 200 equi-spaced bins between their empirical 0.025 and 0.975 quantiles.

The selection of a ridge parameter can be avoided if, instead of using $\hat{f}_U^{\mathrm{Ft}}(t) + \rho$ in $\hat{L}$, we employ $\tilde{f}_U^{\mathrm{Ft}}(t) = \hat{f}_U^{\mathrm{Ft}}(t) I(t \in A) + \hat{f}_P^{\mathrm{Ft}}(t) I(t \notin A)$, where $A$ denotes the largest interval around 0 in which $\hat{f}_U^{\mathrm{Ft}}(t)$ is nonincreasing to the left of 0 and nonincreasing to the right (excepting fluctuations very close to 0), and $\hat{f}_P^{\mathrm{Ft}}(t)$ is a parametric function estimated from the observations and defined by $\hat{f}_P^{\mathrm{Ft}}(t) = (1 + A_U t^2)^{-B_U}$, with $A_U$ and $B_U$ chosen so as to match the empirical second and fourth moments of the error with those of $\hat{f}_P$. In the event that these moments are negative, we set $B_U = 1$ and take $A_U$ equal to half the empirical variance of the error, which corresponds to $\hat{f}_P$ being a Laplace density. This method gives very good results in practice, sometimes even better than in the case of known error density. It is designed specifically for the comparatively small samples that typically arise in errors-in-variables regression with repeated measurements. The small sample sizes there are typically a consequence of the relatively high cost, in terms of time, effort or money, of making several observations of the same $X$, compared with making the same number of observations of different $X$'s.

In our simulations we consider samples of sizes $n = 50$, 100 and 250, and fix the number of replications, $N_j$, at 2 or 4. In each case we generate 500 datasets, for each of which we calculate an estimate of the target curve by using the bandwidth $\hat{h}_{\text{PI}}$ (density case) or the bandwidth $\hat{h}_{\text{CV}}$ (regression case). We take $K^{\text{Ft}} = (1-t^2)^3 I(t \in [-1,1])$; this kernel is commonly used in deconvolution problems. To evaluate performance, we calculate the integrated squared error (ISE) distance of each estimate, where $\text{ISE} = \int_I (\hat{m} - m)^2$, with $m = f_X$ or $m = g$, and where $I$ is the whole real line (density case) or $I = [0,1]$ (regression case). In the graphs we present the three estimates that resulted in the first, second and third quartiles of the 500 calculated ISEs, and we denote them by, respectively, $q_1$, $q_2$ and $q_3$. We report only part of the simulations, although our conclusions are similar for the other, nonreported results.

Table 1 illustrates the effect of increasing the number of replications by comparing the median and the inter-quartile range (IQR) of the calculated ISEs, for $N_j = 2$ and $N_j = 4$, obtained from 500 samples from density (i) contaminated by Laplace or normal errors when $M = 200$ or $M = 500$. These and related results indicate better performance when $N_j = 2$ than when $N_j = 4$. As suggested in the introduction, for the same total number of observations, $M$, it is more advantageous to have a large number of intrinsically different observations, $n$, than a large number of replications, $N_j$.

In Figure 1 we show the quartile curves obtained for 500 samples from density (iii) contaminated by Laplace error when $n = 250$, with $N_j = 2$, together with boxplots of the calculated ISEs for $n = 50$, 100 and 250 in the known and unknown error cases. The results show that, as in the error-free case, it is difficult to recover all the modes of this density. They also illustrate the fact that knowing the error density brings only minor improvements, which we also observed in our non reported simulated results. In some of the non-reported cases, the results were even better for $\hat{f}_X$ than for $\tilde{f}_X$.

Figure 2 shows the quartile curves obtained from 500 samples in the case of regression function (i) for $n = 100$ and $N_j = 2$, when the error $U$ is normal and $X \sim U[0,1]$. We also show boxplots of the 500 calculated ISEs in the case of Laplace and normal error $U$ and $n = 100$ or 250, using $\hat{g}$ (unknown error) or $\tilde{g}$ (known error). We see that the estimated curves are quite good and the results are slightly better when the error density is known.

TABLE 1
*Values of* median $\times$ 100 (IQR $\times$ 100) *of the ISE for density* (i), *when* $M = 200$ *or* $M = 500$, *with* $N_j = 2$ *or* $N_j = 4$ *and Laplace or normal errors*

| $(N_j, M)$   | (2, 200)    | (4, 200)    | (2, 500)    | (4, 500)    |
|--------------|-------------|-------------|-------------|-------------|
| $U \sim$ Lap | 1.41 (0.94) | 1.56 (0.98) | 0.89 (0.51) | 0.96 (0.58) |
| $U \sim$ Norm| 2.09 (1.33) | 2.31 (1.43) | 1.42 (0.92) | 1.55 (1.02) |

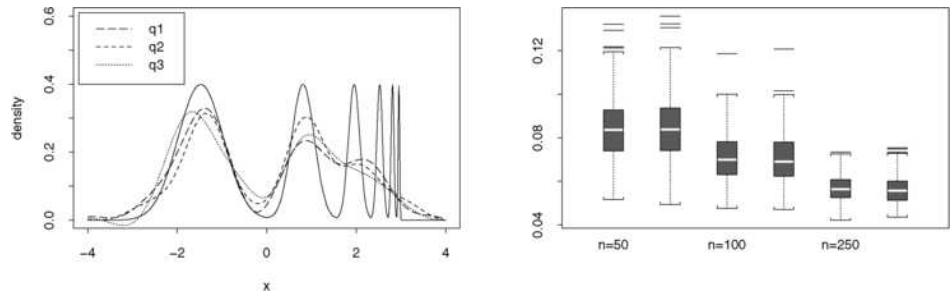

FIG. 1. *Quartile curves of 500 estimates $\hat{f}_X$ of density* (iii) *in the Laplace error case, for $N_j = 2$ and $n = 250$* (*left panel*), *together with boxplots* (*right panel*) *of the 500 calculated ISEs when $n = 50$, 100 or $n = 250$. In each group of two boxplots, the first is for $\hat{f}_X$ and the second for $\tilde{f}_X$.*

Finally, Figure 3 shows the quartile curves in the case of regression curve (iii), when the error $U$ is Laplace, $X \sim U[0,1]$, $N_j = 2$ and $n = 100$ or 250. In this case, too, we see that the results are quite good and improve as sample size increases.

4.2. *Real-data examples.* We apply our methods to two medical examples. The first dataset, described by Bland and Altman (1986), was collected to compare two methods for measuring peak expiratory flow rate (PEFR). Two replicated measurements of PEFR were made on 17 individuals, using each of two different methods: a Wright peak flow meter and a mini Wright meter. As described by Bland and Altman (1986), when evaluating a new method for measuring a clinical quantity, usually the true values remain unknown and a common practice is to compare the new method with the established method, rather than with the true quantities. The goal is thus to check whether the mini meter and the Wright meter are in agreement.

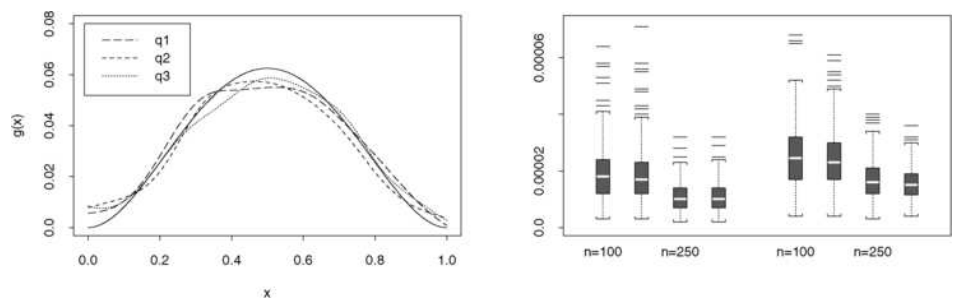

FIG. 2. *Quartile curves of 500 estimates $\hat{g}$ of the regression function* (i) *in the normal error case for $N_j = 2$, $n = 100$ and $X \sim U[0,1]$* (*left panel*); *and boxplots of 500 ISEs for the same regression curve in the case of Laplace error* (*first group of four*) *or normal error* (*last group of four*), *for $n = 100$ or 250* (*right panel*). *In each group of two boxplots, the first is for $\hat{g}$ and the second is for $\tilde{g}$.*

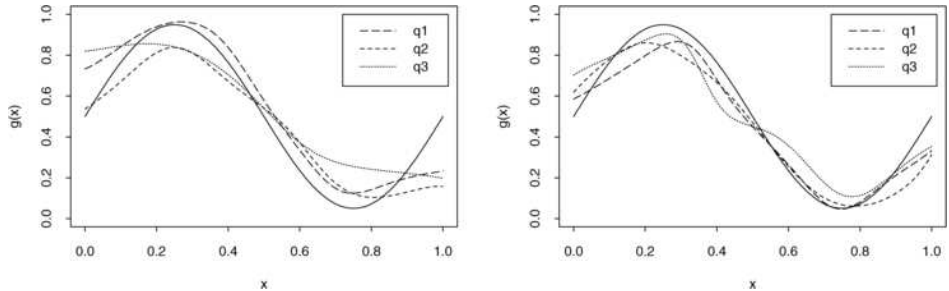

FIG. 3. *Quartile curves of 500 estimates $\hat{g}$ of the regression function* (iii) *in the Laplace error case for $N_j = 2$, $X \sim U[0,1]$ and $n = 100$* (*left panel*) *or $n = 250$* (*right panel*).

To this end, we define $X_i$ as the average of all possible readings on the mini meter for individual $i$, and define $Y_i$ similarly for the "regular" Wright meter. The latter gives more stable (less variable) readings than the mini meter, and, therefore, for each individual $i$, we set $Y_i$ equal to the average of the two Wright readings. Since readings from the mini meter are more variable, then there we need to incorporate measurement errors. For $j = 1, 2$, we take $W_{ij}$ to be the $j$th replicated mini Wright measurement.

The regression estimate (dashed line) is depicted in the left panel of Figure 4, together with the Nadaraya–Watson estimate of $g$ (dotted line) that uses the original data (and hence, ignores the error $U$) and the data $(W_{ij}, Y_i)$. The unusual shape of the dashed line, deviant from a straight line, suggests that the two PERF measurement methods might not be in good agreement and that further investigation should be carried out. Bland and Altman (1986) note that a standard parametric analysis of these data, not taking the noise into account, indicates agreement between the two methods. Analogously, the dotted line shows that ignoring the measurement error results in an estimate that oversmoothes the data, and which lies closer to, although still far from, a straight line. For example, the steeper climb of the dashed line in the upper right-hand part of the graph, and the flatter nature of that line after the climb (compared to the dotted line), add weight to the argument that the results from the mini Wright meter, in the 600–700 range, represent stochastic variation of the relatively constant measurements obtained using the Wright peak flow meter.

The second dataset concerns two replicated measurements derived from CAT scans of the heads of 50 psychiatric patients. More precisely, the ventricule-brain ratio (VBR) was measured twice for each patient, using a hand-held planimeter. See Turner, Toone and Brett-Jones (1986) and Dunn (2004). The logarithm of the VBR can be described by model (2.1), and for the $i$th patient we set $W_{ij} = \log(\text{VBR}_{ij})$, $j = 1, 2$, where $\text{VBR}_{ij}$ denotes the $j$th contaminated replication of the measurement of VBR for patient $i$. The density estimate of the noncontaminated log VBR is plotted as the dashed line in Figure 4, and represents a smooth and symmetric density. We also show, in the dotted line, the kernel density estimate that ignores measurement error. The second estimate is essentially a smoothed version of the density shown by the dashed line, modulo Gibbs-phenomenon wiggles in the tails of the latter.

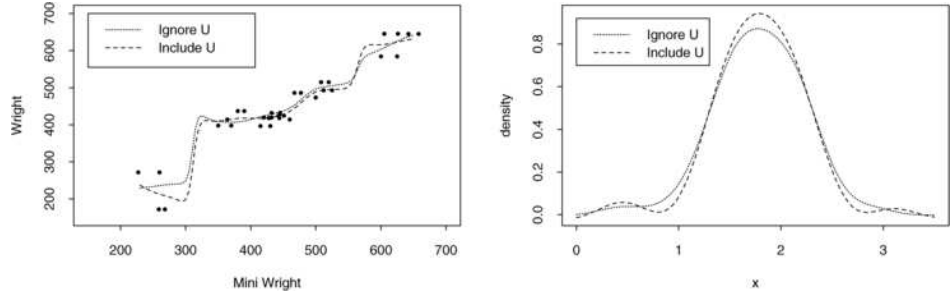

Fig. 4. *Regression estimate for the PEFR data (left panel) and density estimate for the CAT data (right panel).*

**5. Outlines of technical arguments.** Details of proofs, and a derivation of Theorem 3.3, are given by Delaigle, Hall and Meister (2006). Without loss of generality, $c = 1$ in (3.2).

5.1. *Outline proof of Theorem* 3.1. Put $\psi = f_U^{\text{Ft}}$, $\phi = \psi^2$ and

$$\Delta(t) = \frac{1}{N} \sum_{j=1}^n \sum_{(k_1, k_2) \in \mathcal{S}_j} [\cos\{t(W_{jk_1} - W_{jk_2})\} - \phi(t)].$$

In this notation,

$$(5.1) \qquad (\hat{f}_U^{\text{Ft}} + \rho)^{-1} = \psi^{-1} I(\psi > \rho) + \sum_{\ell=1}^k c_\ell \psi^{-2\ell-1} \Delta^\ell + \chi_1 + \chi_2,$$

where the constants $c_\ell$ are derived from binomial coefficients, $|\chi_1| \leq \rho^{-1} I(\psi \leq \rho)$,

$$|\chi_2| \leq \text{const.} \bigg\{ \frac{\rho}{\psi + \rho}(\psi^{-3}|\Delta| + \psi^{-(2k+1)}|\Delta|^k) + \psi^{-(2k+3)}|\Delta|^{k+1} \\ + \rho\psi^{-2} I(\psi > \rho) + \rho^{-1} I\bigg(|\Delta| > \frac{1}{2}\phi\bigg) \bigg\},$$

and "const.," here and below, denotes a generic positive constant depending only on $k$, $f_U$ and the parameters $\alpha$ and $C_2$ of $\mathcal{F}(\beta, C_2)$.

Result (5.1) implies that

$$\text{(5.2)} \quad \hat{f}_X(x) - \tilde{f}_X(x) = \sum_{\ell=1}^{k} c_\ell \delta_{1\ell}(x) + \delta_{01}(x) + \delta_{02}(x) - \delta_2(x),$$

where, for $\ell = 1, 2$ in the case of $\delta_{0\ell}$, and $1 \leq \ell \leq k$ for $\delta_{1\ell}$,

$$\delta_{0\ell}(x) = \frac{1}{2\pi} \int e^{-itx} \hat{f}_W^{\text{Ft}}(t) \chi_\ell(t) K^{\text{Ft}}(ht) \, dt,$$

$$\delta_{1\ell}(x) = \frac{1}{2\pi} \int e^{-itx} \hat{f}_W^{\text{Ft}}(t) \psi(t)^{-2\ell-1} \Delta(t)^\ell K^{\text{Ft}}(ht) \, dt,$$

$$\delta_2(x) = \frac{1}{2\pi} \int e^{-itx} \hat{f}_W^{\text{Ft}}(t) \psi(t)^{-1} K^{\text{Ft}}(ht) I\{\psi(t) \leq \rho\} \, dt$$

and $\hat{f}_W^{\text{Ft}}(t) = M^{-1} \sum_j \sum_k e^{itW_{jk}}$.

It can be proved that, for a constant $n_0 \geq 1$, the functions $\delta_{01}$ and $\delta_2$ vanish identically whenever $n \geq n_0$. Therefore, assuming $n \geq n_0$, we deduce from (5.2) that

$$\hat{f}_X(x) - \tilde{f}_X(x) = \sum_{\ell=1}^{k} c_\ell \delta_{1\ell}(x) + \delta_{02}(x).$$

This formula and the fact that $\hat{f}_W^{\text{Ft}} = \psi f_X^{\text{Ft}} + \Delta_1$, where $\Delta_1 = \hat{f}_W^{\text{Ft}} - E(\hat{f}_W^{\text{Ft}})$, imply that

$$\text{(5.3)} \quad \sup_{-\infty < x < \infty} E\{\hat{f}_X(x) - \tilde{f}_X(x)\}^2$$

$$\leq \text{const.} \left[ \max_{r=2,3} \max_{1 \leq \ell \leq k} \sup_{-\infty < x < \infty} E\{\delta_{r\ell}(x)^2\} + \sup_{-\infty < x < \infty} E\{\delta_{02}(x)^2\} \right],$$

where

$$\delta_{2\ell}(x) = \frac{1}{2\pi} \int e^{-itx} f_X^{\text{Ft}}(t) \psi(t)^{-2\ell} \Delta(t)^\ell K^{\text{Ft}}(ht) \, dt,$$

$$\delta_{3\ell}(x) = \frac{1}{2\pi} \int e^{-itx} \psi(t)^{-2\ell-1} \Delta_1(t) \Delta(t)^\ell K^{\text{Ft}}(ht) \, dt.$$

Lengthy arguments can be used to show that

$$\max_{r=2,3} \max_{1 \leq \ell \leq k} \sup_{-\infty < x < \infty} E\{\delta_{r\ell}(x)^2\}$$

$$\leq \text{const.}[n^{-1}\{h^{\beta - 2\alpha - 1} + h^{2(\beta - 2\alpha) - 1} + (\log n)^2\}$$

$$+ n^{-2}(h^{2(\beta - 4\alpha) - 2} + h^{-6\alpha - 1}) + n^{-k} h^{2\beta - 4k\alpha - 2}$$

$$+ n^{-(k+1)} h^{-2(2k+1)\alpha - 1}]$$

$$= O(p_n)$$

and $\sup_x E\{\delta_{02}(x)^2\} = O(p_n)$. Together, these bounds and (5.3) imply (3.3).

5.2. *Outline proof of Theorem* 3.3. For brevity we derive only the second part of (3.6). Since $|\{\hat{f}_U^{\mathrm{Ft}}(t) + \rho\}^{-1} - \hat{f}_U^{\mathrm{Ft}}(t)^{-1}| \leq \rho/\hat{f}_U^{\mathrm{Ft}}(t)^2$, then

$$(5.4) \qquad |\widehat{L}(u) - \widehat{L}^0(u)| \leq \frac{\rho h}{2\pi} \int \hat{f}_U^{\mathrm{Ft}}(t)^{-2} K^{\mathrm{Ft}}(ht)\, dt,$$

where $\widehat{L}^0$ denotes the version of $\widehat{L}$ constructed with $\rho = 0$. With probability $\pi_n$, say, equal to $1 - O(n^{-B})$ for each $B > 0$, $\frac{1}{2} f_U^{\mathrm{Ft}}(t)^2 \leq \hat{f}_U^{\mathrm{Ft}}(t)^2$ for all $t$ such that the integrand at (5.4) does not vanish. Therefore, with probability at least $\pi_n$,

$$\sup_{-\infty < u < \infty} |\widehat{L}(u) - \widehat{L}^0(u)| \leq \frac{C_1^2 \rho h s}{\pi} \int_{-1/h}^{1/h} (1 + |t|)^{2\alpha}\, dt \leq C_3 \rho h^{-2\alpha},$$

where $s = \sup |K^{\mathrm{Ft}}|$ and $C_3 > 0$. Hence, with probability at least $\pi_n$,

$$\sup_{-\infty < x < \infty} |\hat{f}_X(x) - \hat{f}_X^0(x)| \leq C_3 n^{-1},$$

which leads to the second part of (3.6).

A. Delaigle
Department of Mathematics
University of Bristol
Bristol B98 4JS
United Kingdom
E-mail: Aurore.Delaigle@bristol.ac.uk

P. Hall
Department of Mathematics
and Statistics
University of Melbourne
Victoria 3010
Australia

A. Meister
Institut für Stochastik
und Anwendungen
Universität Stuttgart
D-70569 Stuttgart
Germany